\documentclass[12pt,oneside,reqno]{amsproc}

\usepackage[cp1251]{inputenc}
\usepackage[russian]{babel}
\usepackage{amssymb,amsmath,amsthm}

\numberwithin{equation}{section}
\newtheorem{lem}{Лемма}[section]
\newtheorem{tm}{Теорема}[section]
\newtheorem{sek}{Следствие}[section]

\newtheorem{rim}{Замечание}[section]
\newtheorem{dif}{Определение}[section]

\newcommand{\Var}{\mathop{\rm Var}\nolimits}

\title[]{О конструкции и некоторых свойствах самоподобных функций в пространствах $L_p[0,1]$}
\author{И.~А.~Шейпак\footnote{%
Работа выполнена при поддержке грантов РФФИ \No\,04-01-00712 и
поддержки ведущих научных школ \No\,НШ-5247.2006.1}%
}

\begin{document}
\noindent УДК 517.518
\begin{abstract}
В статье приводится конструкция аффинно-самоподобных функций. В терминах параметров преобразований
самоподобия дается условие принадлежности этих функций классам $L_p[0,1]$, а также пространству
$C[0,1]$. Изучаются некоторый свойства этих функций (монотонность, ограниченность вариации).
Устанавливается связь самоподобных функций с самоподобными мерами.
\end{abstract}
\maketitle

\section{Введение} Самоподобные (фрактальные) объекты (множества, меры, функции) нашли широкое применение
в различных областях математики. Одним из способов построения таких объектов является применение
теоремы о неподвижной точке сжимающего отображения в полном метрическом пространстве. Достаточно
общую конструкцию самоподобных множеств и мер, опирающуюся на такой подход, а также некоторые их
приложения можно найти, например, в~\cite{Hu}.

Одним из важных приложений самоподобных непрерывных функций является теория сжатия компьютерных
изображений~(см. \cite{Bar}). При этом, в основном, решаются вопросы интерполяции изображаемых
объектов или их границ непрерывными самоподобными функциями, значения которых заданы в конечном
наборе точек. В конструкции этих функций параметры отображений, определяющих самоподобие, могут
быть вычислены через узлы интерполяции.

Кроме того, самоподобные меры находят применение и в спектральной теории операторов~(\cite{SV}) и в
гармоническом анализе~(\cite{Str},~\cite{Zig}). Заметим, что в работе~\cite{Zig} терминология
"<самоподобные меры"> не употреблялась, но конструкция мер, коэффициенты Фурье-Стильтьеса которых
не стремятся к нулю, тем не менее использует идею самоподобия. Различные задачи, связанные с
изучением свойств вероятностных самоподобных мер, естественным образом возникают в теории
фрактальных кривых, изучение которых началось с уже ставших классическими работы
\cite{deRam1}--\cite{deRam3}. Впоследствии непрерывные фрактальные кривые нашли применение в теории
функций, теории вероятностей, эргодической теории и т.д.~\cite{CavDah}--\cite{Nik}.

Большое развитие  теория фрактальных кривых получила после того, как обнаружилась её связь с
теорией всплесков (вейвлетов) и масштабирующих функций (см., например, \cite{Daub1}). В частности,
определённый интерес представляет исследование гладкости решений масштабирующих уравнений в
различных функциональных пространствах. Например, в работе \cite{Daub2} были получены оценки сверху
на показатель Г\"ельдера этих решений в пространстве $C[0,1]$. В~\cite{Prot1} были получены
критерии таких свойств решений масштабирующих уравнений, как абсолютная непрерывность, сингулярная
непрерывность, ограниченность вариации.

Расширением понятий самоподобных мер и непрерывных функций являются функции из пространств $L_p$. В
связи с этим необходимо упомянуть о введ\"енных в работе~\cite{Prot2} суммируемых фрактальных
кривых. В этой работе были получены критерии существования фрактальной кривой и принадлежности е\"е
классам $L_p$ в терминах спектральных $p$-радиусов $\rho_p$~\cite{LauWang}. Кроме того, там же были
выведены формулы для показателей гладкости в различных функциональных пространствах.

Некоторые частные случаи квадратично-суммируемых самоподобных функций были рассмотрены
в~\cite{VlSh}. Свойства самих функций не изучались. Основное отличие конструкции самоподобных
функций из $L_p[0,1]$ от самоподобных непрерывных функций заключается в том, что параметры
самоподобия нельзя определить через значения функции в некотором наборе точек.

Упомянутые работы (\cite{SV}, \cite{VlSh}), а также некоторые результаты работы~\cite{KacKr}
указывают на тесную связь спектра некоторых граничных задач с такими свойствами самоподобных
функций как монотонность, непрерывность, абсолютная непрерывность, ограниченность вариации. В связи
с этим изучение указанных и других свойств самоподобных функций различных классов представляет
интерес для многих математических направлений.

Целью данной работы является построение самоподобных функций в $L_p[0,1]$ и исследование некоторых
их свойств (непрерывность, монотонность, ограниченность вариации). В основе конструирования таких
функций лежит теорема о неподвижной точке сжимающего отображения в полном метрическом пространстве.
Предлагаемый подход, с одной стороны является более общим, по сравнению с построением фрактальных
кривых (отрезки разбиения могут иметь различную длину). С другой стороны, фрактальная кривая
принадлежит $\mathbb R^n$. В случае, когда самоподобная функция является непрерывной, эту
конструкцию можно переформулировать в терминах узлов интерполяции (см.~\cite{Bar}). В случае, когда
самоподобная функция является неубывающей, она естественным образом порождает самоподобную меру.
Самоподобной функции ограниченной вариации соответствует самоподобный заряд.

\section{Самоподобные функции в пространстве \(L_p[0,1]\)}\label{pt:2}
\subsection{Операторы подобия в пространстве \(L_p[0,1]\)}
Пусть фиксировано натуральное число \(n>1\), и пусть вещественные числа
\(a_k>0\), \(c_k>0\), \(d_k\) и \(\beta_k\), где \(k=1,\ldots,n\), таковы, что
\[
\sum\limits_{k=1}^n a_k=1.
\]
Данному набору чисел и произвольному $x\in[0,1]$ можно поставить в соответствие
непрерывный нелинейный оператор \(G:L_p[0,1]\to L_p[0,1]\) вида
\begin{equation}\label{eq:auxto}
G(f)=\sum\limits_{k=1}^n\left\{\beta_k\cdot\chi_{(\alpha_k,\alpha_{k+1})}+c_k\cdot
x +d_k\cdot G_k(f)\right\},
\end{equation}
где использованы следующие обозначения:
\begin{enumerate}
\item через $\alpha_k$, где $k=1,2,\ldots,n+1$, обозначены числа
$\alpha_1=0$ и $\alpha_k=\sum_{l=1}^{k-1}a_l$, где $k=2,\ldots,n+1$;
\item через $\chi_{(\zeta,\xi)}$ обозначена
характеристическая функция интервала $(\zeta,\xi)$, рассматриваемая как элемент пространства
$L_p[0,1]$; \item через $G_k$, где $k=1,\ldots,n$, обозначены непрерывные линейные операторы в
пространстве $L_p[0,1]$, действующие на функцию $f$ согласно правилу
\begin{equation}\label{eq3:2.1}
G_k(f)(t)=f(a_kt+\alpha_k),
\qquad k=1,\ldots,n.
\end{equation}
\end{enumerate}

Операторы \(G\) вида~\eqref{eq:auxto} будут называться \emph{операторами подобия}.

Этот оператор на непрерывные функции действует  следующим образом:
\begin{equation}\label{eq:nepr_otobr}
G(f)(t)=c_k\cdot t+d_k\cdot f(t)+\beta_k,
\end{equation}
где $x\in[\alpha_{k},\alpha_{k+1}]$, $t\in[0,1]$, а переменные $x$ и $t$ связаны соотношением
$t=a_k\cdot x+\alpha_k$, $k=1,2,\ldots,n$.

Это отображение является частью б\'ольшего отображения $A$, действующего на
график непрерывной функции. Чтобы описать его, введем операторы $A_k$:
$$
A_k\left(\begin{array}{c}t\\f(t)\end{array}\right)=
\left(\begin{array}{cc} a_k & 0\\ c_k &
d_k\end{array}\right)
\left(\begin{array}{c}t\\f(t)\end{array}\right)+
\left(\begin{array}{c}\alpha_k\\ \beta_k\end{array}\right).
$$
Отображение $A$ определено так
$$
A\left(\begin{array}{c}t\\f(t)\end{array}\right)=\bigcup_{k=1}^n A_k
\left(\begin{array}{c}t\\f(t)\end{array}\right).
$$

Среди всех операторов подобия для дальнейшего изложения наибольший интерес представляют сжимающие
операторы.

\begin{lem}\label{lem2:1}
Оператор подобия $G$ является сжимающим в $L_p[0,1]$ в том и только том случае, когда справедливы
неравенства
\begin{gather}\label{eq:szim}
\sum\limits_{k=1}^n a_k\,|d_k|^p<1 \quad (1\leqslant p<+\infty),\\\label{eq:szim2}
\max_{1\leqslant k\leqslant n}|d_k|<1 \quad (p=+\infty).
\end{gather}
\end{lem}
\begin{proof}
При $p<\infty$ утверждение доказываемой леммы следует из того факта, что при любых
\(f_1,f_2\in L_p[0,1]\) справедливы  соотношения
\begin{multline*}
\|G(f_1)-G(f_2)\|_{L_p[0,1]}^p=\int\limits_0^1|G(f_1)-G(f_2)|^p\,dx=\\\sum\limits_{k=1}^n\left(|d_k|^p\,\int\limits_{\alpha_k}^{\alpha_{k+1}}
|G_k(f_1)-G_k(f_2)|^p\,dx\right)=\left(\sum\limits_{k=1}^n a_k\,|d_k|^p
\right)\,\int\limits_0^1|f_1-f_2|^p\,dx=\\=\left(\sum\limits_{k=1}^n a_k\,
|d_k|^p\right)\,\|f_1-f_2\|_{L_2[0,1]}^p.
\end{multline*}
При $p=+\infty$ справедливо  неравенство
$$
\|f_1-f_2\|_{L_\infty[0,1]}\leqslant \max_{k}|d_k|\|f_1-f_2\|_{L_\infty[0,1]}.
$$
\end{proof}

Из леммы~\ref{lem2:1} и принципа сжимающих отображений немедленно вытекает справедливость
следующего утверждения.
\begin{tm}\label{sek2:1}
Если справедливо неравенство~\eqref{eq:szim} ($1\leqslant p<+\infty$), или~\eqref{eq:szim2}
($p=+\infty$), то существует и единственна функция $f\in L_p[0,1]$, удовлетворяющая уравнению
$G(f)=f$.
\end{tm}

\begin{dif} Функции, заданные условием $G(f)=f$ при некотором сжимающем
операторе подобия $G$ будем называть аффинно-самоподобными или просто самоподобными. Числа
 $\{a_k\}$, $\{c_k\}$, $\{d_k\}$ и $\{\beta_k\}$, $k=1,2,\ldots,n$ будем называть параметрами
 самоподобия.
\end{dif}

В дальнейшем всегда будет предполагаться, что неравенство~\eqref{eq:szim} или~\eqref{eq:szim2}
выполнено.

\begin{rim}
Набор чисел $\{\alpha_k\}$, как нетрудно видеть, задает разбиение отрезка $[0,1]$:
$0=\alpha_1<\alpha_2<\ldots<\alpha_n<1=\alpha_{n+1}$. Числа $a_k$ отвечают за горизонтальное сжатие
функции $f$, а числа $d_k$ за вертикальное сжатие.
\end{rim}

\begin{rim}
Параметры $c_k$ не влияют на свойство сжимаемости оператора $G$.
\end{rim}

\begin{rim}
Для фрактальных кривых, заданных семейством аффинных (не обязательно сжимающих) операторов условию
\eqref{eq:szim} при $p\in [1,+\infty)$ соответствует условие $\rho_p<1$ (см., например, ~\cite{Prot2}), где $\rho_p$
--- $p$-радиус этих операторов.
\end{rim}

Одной и той же функции могут соответствовать разные параметры самоподобия. Например, функция
$f(x)=x$ на отрезке $[0,1]$ может быть задана как набором $n=2$, $a_1=a_2=1/2$, $c_1=c_2=1/2$,
$d_1=d_2=0$, $\beta_1=0$, $\beta_2=1/2$, так и  набором $n=3$, $a_1=a_2=a_3=1/3$,
$c_1=c_2=c_3=1/3$, $d_1=d_2=d_3=0$, $\beta_1=0$, $\beta_2=1/3$, $\beta_3=1/3$.

\subsection{Неравенства, оценивающие нормы самоподобных функций $f\in L_p[0,1]$ через параметры самоподобия}
На пространстве $n$-мерных векторов $x=(x_1,x_2,\ldots, x_n)$ введем весовую норму для пары
векторов $x=(x_1,x_2,\ldots, x_n)$ и $y=(y_1,y_2,\ldots, y_n)$ и фиксированного $s\in[1,+\infty)$
$$
\|\{x,y\}\|_{s,a}=\left(\sum_{k=1}^n(|x_k|+|y_k|)^sa_k\right)^{1/s},
$$
где коэффициенты $a_k$ удовлетворяют условиям $0<a_k<1$, $k=1,2,\ldots,n$, $\sum_{k=1}^n a_k=1$.

Кроме того, для облегчения выкладок введем обозначение
$$
r_p=\sum_{k=1}^n a_k|d_k|^p.
$$
Легко заметить, что из условия $r_p<1$ для  любого $1\leqslant s<p$ также вытекает условие $r_s<1$
и, более того, $r_s\leqslant (r_p)^{\frac{s}{p}}$.

При $p=+\infty$ введенные выше величины не зависят от чисел $\{a_k\}$ и,  соответственно,
определяются так:
$$
\|\{x,y\}\|_{\infty}=\max_{1\leqslant k\leqslant n}\{|x_k|+|y_k|\},\quad
r_{\infty}=\max_{1\leqslant k\leqslant n}\{|d_k|\}.
$$

Справедлива следующая теорема.
\begin{tm}\label{tm:norm}
Для самоподобной функции $f\in L_p[0,1]$ ($p\in \mathbb N$) с параметрами самоподобия $a_k$, $c_k$,
$d_k$, $\beta_k$  выполнено неравенство
\begin{equation}\label{eq:norm[p]}
\|f\|_{L_p[0,1]}\leqslant \dfrac{\sum_{s=1}^p \|\{c,\beta\}\|_{s,a}}{\left(\prod_{s=1}^p(1-r_s)\right)^{1/p}},
\end{equation}
где $c=(c_1,c_2,\ldots,c_n)$, $\beta=(\beta_1,\beta_2,\ldots,\beta_n)$.
\end{tm}

\begin{proof} Доказательство проведем методом математической индукции по $p$.

Действительно, при $p=1$ справедливы следующие выкладки
\begin{multline*}
\|f\|_{L_1[0,1]}=\int_0^1|f(x)|\,dx=\sum_{k=1}^n\int_{\alpha_k}^{\alpha_{k+1}}|f(x)|\,dx=\\
=\sum_{k=1}^n\int_0^1|c_k\cdot t+\beta_k+d_k\cdot f(t)|a_k\,dt\leqslant
\sum_{k=1}^n (|c_k|+|\beta_k|)a_k+\left(\sum_{k=1}^n a_k|d_k|\right)\|f\|_{L_1[0,1]},
\end{multline*}
из которых следует справедливость неравенства~\eqref{eq:norm[p]} при $p=1$.

Допустим, что неравенство~\eqref{eq:norm[p]} верно при $j=1,2,\ldots, p-1$. Если $f\in L_p[0,1]$,
то, очевидно, $f\in L_s[0,1]$ $\forall\; s\in[1,p)$. Тогда справедливы следующие преобразования
\begin{multline*}
\|f\|^p_{L_p[0,1]}=\sum_{k=1}^n\int_{\alpha_k}^{\alpha_{k+1}}|f(x)|^p\,dx=
\sum_{k=1}^n\int_0^1|c_k\cdot t+\beta_k+d_k\cdot f(t)|^pa_k\,dt\leqslant\\ \leqslant
\sum_{k=1}^n\sum_{j=0}^{p-1} C_p^j (|c_k|+|\beta_k|)^{p-j}|d_k|^ja_k\|f\|^j_{L_j[0,1]}+
\left(\sum_{k=1}^n a_k|d_k|^p\right)\|f\|^p_{L_p[0,1]},
\end{multline*}
из которых вытекает неравенство
\begin{equation}\label{eq:ner1}
(1-r_p)\|f\|^p_{L_p[0,1]}\leqslant
\sum_{k=1}(|c_k|+|d_k|)^{p}a_k+\sum_{k=1}^n\sum_{j=1}^{p-1} C_p^j
(|c_k|+|d_k|)^{p-j}|d_k|^ja_k\|f\|^j_{L_j[0,1]}.
\end{equation}
Учитывая предположение индукции при $j=1,2,\ldots, p-1$:
$$
\|f\|^j_{L_j[0,1]}\leqslant\dfrac{\left(\sum_{s=1}^j\|\{c,\beta\}\|_{s,a}\right)^j}{\prod_{s=1}^j(1-r_s)}
\leqslant
\dfrac{\left(\sum_{s=1}^{p-1} \|\{c,\beta\}\|_{s,a}\right)^j}{\prod_{s=1}^{p-1}(1-r_s)},
$$
а также применяя неравенство Г\"ельдера к сумме:
\begin{multline*}
\sum_{k=1}^n (|c_k|+|d_k|)^{p-j}|d_k|^ja_k=\sum_{k=1}^n
(|c_k|+|d_k|)^{p-j}a_k^{\frac{p-j}{p}}a_k^{\frac{p}{j}}|d_k|^j\leqslant\\
\leqslant \left(\sum_{k=1}^n (|c_k|+|d_k|)^{p}a_k\right)^{\frac{p-j}{p}}
\left(\sum_{k=1}^n a_k|d_k|^p\right)^{\frac{j}{p}}=\|\{c,\beta\}\|^{p-j}_{p,a}r_p^{\frac{j}{p}}\leqslant
\|\{c,\beta\}\|^{p-j}_{p,a},
\end{multline*}
неравенство~\eqref{eq:ner1} преобразуем к следующему виду:
$$
(1-r_p)\|f\|^p_{L_p[0,1]}\leqslant \|\{c,\beta\}\|^{p}_{p,a}+\sum_{j=1}^{p-1}
C_p^j\|\{c,\beta\}\|^{p-j}_{p,a}
\dfrac{\left(\sum_{s=1}^{p-1}\|\{c,\beta\}\|_{s,a}\right)^j}{\prod_{s=1}^{p-1}(1-r_s)}.
$$
Из последнего неравенства тем более следует неравенство
$$
(1-r_p)\|f\|^p_{L_p[0,1]}\leqslant
\dfrac{\left(\sum_{s=1}^{p}\|\{c,\beta\}\|_{s,a}\right)^p}{\prod_{s=1}^{p-1}(1-r_s)},
$$
что и завершает доказательство теоремы.
\end{proof}

\begin{rim}
В пространстве $L_\infty[0,1]$ неравенство для нормы $f$, аналогичное~\eqref{eq:norm[p]} имеет вид:
$$
\|f\|_{L_\infty[0,1]}\leqslant\dfrac{\max_{k}\{|c_k|+|\beta_k|\}}{1-\max_{k}|d_k|}=
\dfrac{\|\{c,\beta\}\|_\infty}{1-r_\infty}.
$$
\end{rim}

Докажем аналогичные оценки для дробного $p$. Для этого нецелое число $p> 1$ представим в виде
$p=[p]+\{p\}$, где $[p]$ --- обозначает целую часть числа $p$, а $\{p\}$ --- дробную часть числа
$p$, $\{p\}\ne 0$.

\begin{tm}\label{tm:norm_frac}
Для самоподобной функции $f\in L_p[0,1]$, $\{p\}\ne 0$ справедлива оценка
$$
\|f\|_{L_p[0,1]}\leqslant
C\cdot
\dfrac{\left(\|\{c,\beta\}\|_{[p],a}^{[p]}+
\sum_{s=1}^{[p]}\|\{c,\beta\}\|_{s,a}\right)^\frac{[p]}{p}}{\left((1-r_p)\prod_{s=1}^{[p]}(1-r_s)\right)^\frac{1}{p}}.
$$
где
$$
C=\left(\max\left\{\max_{k}\{(|c_k|+|\beta_k|)^{\{p\}}\},\max_k\{|d_k|^{\{p\}}\},
\left(\sum_{s=1}^{[p]}\|\{c,\beta\}\|_{s,a}\right)^{\{p\}}\right\}\right)^\frac{1}{p}.
$$

\end{tm}
\begin{proof}
Справедлива следующая цепочка неравенств
\begin{multline}\label{eq:otzenka_p}
\|f\|_{L_{p}[0,1]}\leqslant\sum_{k=1}^{n}\int_0^1|c_k\cdot t+\beta_k+d_k\cdot f(t)|^{[p]}
(|c_k|\cdot t+|\beta_k|+|d_k|\cdot |f(t)|)^{\{p\}}a_k\,dt\leqslant\\
\leqslant
\max_{k}\{(|c_k|+|\beta_k|)^{\{p\}}\}\sum_{k=1}^{n}
\int_0^1|c_k\cdot t+\beta_k+d_k\cdot f(t)|^{[p]}a_k\,dt+\\+
\sum_{k=1}^{n}
\int_0^1|c_k\cdot t+\beta_k+d_k\cdot f(t)|^{[p]}|d_k|^{\{p\}}|f(t)|^{\{p\}}a_k\,dt\leqslant\\
\leqslant
\max_{k}\{(|c_k|+|\beta_k|)^{\{p\}}\}\|f\|_{L_{[p]}[0,1]}^{[p]}+
\sum_{k=1}^{n}\int_0^1(|c_k|+|\beta_k|+|d_k|\cdot |f(t)|)^{[p]}|d_k|^{\{p\}}|f(t)|^{\{p\}}a_k\,dt.
\end{multline}
Преобразуем последнюю сумму в последнем неравенстве.
\begin{multline*}
\sum_{k=1}^{n}\int_0^1(|c_k|+|\beta_k|+|d_k|\cdot
|f(t)|)^{[p]}|d_k|^{\{p\}}|f(t)|^{\{p\}}a_k\,dt=\\=
\sum_{k=1}^{n}\int_0^1\left(\sum_{j=0}^{[p]}
C_{[p]}^j(|c_k|+|\beta_k|)^{[p]-j}|d_k|^j|f(t)|^j\right)|d_k|^{\{p\}}|f(t)|^{\{p\}}a_k\,dt.
\end{multline*}
В последнем выражении слагаемое при $j=[p]$ имеет вид
\begin{equation}\label{eq:rp}
\sum_{k=1}^{n} a_k|d_k|^p\int_0^1 |f(t)|^p\,dt=r_p\|f\|^p_{L_p[0,1]}.
\end{equation}

Поменяв порядок суммирования, оставшиеся слагаемые можно переписать в виде
\begin{equation}\label{eq:poryadok_sum}
\sum_{j=0}^{[p]-1}
C_{[p]}^j\sum_{k=1}^{n}(|c_k|+|\beta_k|)^{[p]-j}|d_k|^j|d_k|^{\{p\}}a_k\int_0^1|f(t)|^{j+\{p\}}\,dt.
\end{equation}
Суммы по $k$ в выражении~\eqref{eq:poryadok_sum} допускает оценку
\begin{multline}\label{eq:summy}
\sum_{k=1}^{n}(|c_k|+|\beta_k|)^{[p]-j}|d_k|^j|d_k|^{\{p\}}a_k\leqslant
\max_k\{|d_k|^{\{p\}}\}\left(\sum_{k=1}^{n}(|c_k|+|\beta_k|)^{[p]}a_k\right)^{\frac{[p]-j}{[p]}}
\left(\sum_{k=1}^{n}a_k|d_k|^{[p]}\right)^{\frac{j}{[p]}}=\\=
\max_k\{|d_k|^{\{p\}}\}\|\{c,\beta\}\|_{[p],a}^{[p]-j}r_{[p]}^{\frac{j}{[p]}}\leqslant
\max_k\{|d_k|^{\{p\}}\}\|\{c,\beta\}\|_{[p],a}^{[p]-j}, \quad j=0,1,\ldots,[p]-1.
\end{multline}
Здесь мы применили неравенство Г\"ельдера с показателями $p_1=\frac{[p]}{[p]-j}$ и
$q_1=\frac{[p]}{j}$.

Интегралы допускают следующую оценку
\begin{equation}\label{eq:integral}
\int_0^1|f(t)|^{j+\{p\}}\,dt\leqslant \left(\int_0^1|f(t)|^{j+1}\right)^{\frac{j+\{p\}}{j+1}}=
\|f\|^{j+\{p\}}_{L_{j+1}[0,1]}\leqslant\|f\|_{L_{[p]}[0,1]}^{j+\{p\}},
\quad j=0,1,\ldots,[p]-1.
\end{equation}
Здесь также применено неравенство Г\"ельдера с показателями $p_1=\frac{j+1}{j+\{p\}}$ и
$q_1=\frac{j+1}{1-\{p\}}$.

Подставив теперь соотношение~\eqref{eq:rp} и неравенства~\eqref{eq:summy}--\eqref{eq:integral} в
неравенство~\eqref{eq:otzenka_p}, получим
$$
(1-r_p)\|f\|^p_{L_p[0,1]}\leqslant\max_{k}\{(|c_k|+|\beta_k|)^{\{p\}}\}\|f\|_{[p]}^{[p]}+
\max_k\{|d_k|^{\{p\}}\}\sum_{j=0}^{[p]-1}C_{[p]}^j
\|\{c,\beta\}\|_{[p],a}^{[p]-j}\|f\|^{j+\{p\}}_{L_{[p]}[0,1]}.
$$

C учетом соотношения~\eqref{eq:norm[p]} получаем, что
\begin{multline*}
(1-r_p)\|f\|^p_{L_p[0,1]}\leqslant\max_{k}\{(|c_k|+|\beta_k|)^{\{p\}}\}
\dfrac{\left(\sum_{s=1}^{[p]}\|\{c,\beta\}\|_{s,a}\right)^{[p]}}{\prod_{s=1}^{[p]}(1-r_s)}+\\+
\max_k\{|d_k|^{\{p\}}\}\sum_{j=0}^{[p]-1}C_{[p]}^j\|\{c,\beta\}\|_{[p],a}^{[p]-j}
\dfrac{\left(\sum_{s=1}^{[p]}\|\{c,\beta\}\|_{s,a}\right)^{j+\{p\}}}
{\left(\prod_{s=1}^{[p]}(1-r_s)\right)^{\frac{j+\{p\}}{[p]}}},
\end{multline*}
что окончательно можно оценить таким образом:
\begin{multline*}
(1-r_p)\|f\|^p_{L_p[0,1]}\leqslant
\max\left\{\max_{k}\{(|c_k|+|\beta_k|)^{\{p\}}\},\max_k\{|d_k|^{\{p\}}\},
\left(\sum_{s=1}^{[p]}\|\{c,\beta\}\|_{s,a}\right)^{\{p\}}\right\}\times\\\times
\dfrac{\left(\|\{c,\beta\}\|_{[p],a}^{[p]}+
\sum_{s=1}^{[p]}\|\{c,\beta\}\|_{s,a}\right)^{[p]}}{\prod_{s=1}^{[p]}(1-r_s)}.
\end{multline*}
\end{proof}

\begin{sek}\label{sek:1}
Для фиксированного $p\in[1,+\infty]$ и любых $R>0$, $\varepsilon>0$ рассмотрим всевозможные
сжимающие отображения в $L_p[0,1]$, заданные параметрами $\{a_k\}$, $\{c_k\}$, $\{d_k\}$,
$\{\beta_k\}$, $k=1,2,\ldots,n$, такими что
\begin{gather*}
\|\{c,\beta\}\|_{p,a}\leqslant R, \qquad \sum_{k=1}^n a_k|d_k|^p\leqslant 1-\varepsilon,\quad 1\leqslant p<\infty\\
\|\{c,\beta\}\|_{\infty}\leqslant R, \qquad\max_{1\leqslant k\leqslant n}|d_k|\leqslant 1-\varepsilon, \quad p=\infty.
\end{gather*}

Тогда все функции, являющиеся неподвижными точками таких сжимающих отображений, образуют
ограниченное множество в $L_p[0,1]$.
\end{sek}

\subsection{Непрерывная зависимость неподвижной точки сжимающего отображения от параметров
$d_k$, $c_k$, $\beta_k$, $k=1,2,\ldots,n$.}

Пусть задано разбиение отрезка $[0,1]$: $0=\alpha_1<\alpha_2<\ldots<\alpha_n<1$. Числа $a_k$
являются длинами отрезков разбиения:
$$
a_k=\alpha_k-\alpha_{k-1}, \quad k=1,2,\ldots,n.
$$
Зададим также два набора чисел $d_k$, $c_k$, $\beta_k$ и $d'_k$, $c'_k$, $\beta'_k$,
$k=1,2,\ldots,n$, удовлетворяющие условиям
$$
r_p=\sum_{k=1}^n a_k|d_k|^p<1,\qquad r'_p=\sum_{k=1}^n a_k|d'_k|^p<1
$$
при некотором фиксированном $p\in[1,+\infty)$. Эти наборы порождают два сжимающих отображения $G$ и
$G'$, неподвижными точками которых соответственно будут функции $f$  и $g$.

При $p=+\infty$ предполагается, что параметры самоподобия удовлетворяют условиям
$$
r_\infty=\max_{1\leqslant k \leqslant n}\{|d_k|\}<1,\qquad r'_\infty=
\max_{1\leqslant k\leqslant n}\{|d_k|\}<1
$$

Утверждение теорем~\ref{tm:norm} и~\ref{tm:norm_frac} означает, что операторы подобия $G$
непрерывно зависят от параметров $d_k$, $c_k$, $\beta_k$, $k=1,2,\ldots,n$. Применяя известный
результат о непрерывной зависимости неподвижной точки сжимающего оператора от параметра (см.,
например,~\cite{Bar}, лемма 2, стр. 111), получаем следующее утверждение.

\begin{tm}
Самоподобная функция, являющаяся неподвижной точкой сжимающего отображения, непрерывно зависит от
параметров самоподобия, а именно, если $c_k\to c'_k$, $d_k\to d'_k$ и  $\beta_k\to \beta'_k$,
$k=1,2,\ldots,n$, то $\|f-g\|_{L_p[0,1]}\to 0$.
\end{tm}

В частности, при $p\in[1,+\infty)$ несложно получить следующие оценки:
$$
\|f-g\|_{L_p[0,1]}\leqslant\dfrac{2^p\|\{c-c',\beta-\beta'\}\|_{p,a}+2^{p-1}\left(\sum_{k=1}^na_k|d_k-d'_k|^p\right)^\frac{1}{p}
(\|g\|_{L_p[0,1]}+\|f\|_{L_p[0,1]})}{2-r_p^\frac{1}{p}-{r'}_p^\frac{1}{p}}.
$$

В  $L_\infty[0,1]$ аналогичные оценки выглядят так:
$$
\|f-g\|_{L_{\infty}[0,1]}\leqslant\dfrac{2\|\{c-c',\beta-\beta'\}\|_{\infty}+
\max_k\{|d_k-d'_k|\}(\|f\|_{L_\infty[0,1]}+\|g\|_{L_\infty[0,1]})}{2-r_\infty-r'_\infty}.
$$

\section{Непрерывные самоподобные функции}
Укажем условия на числа $\{c_k\}_{k=1}^n$, $\{d_k\}_{k=1}^n$ и $\{\beta_k\}_{k=1}^n$, при которых
оператор подобия $G$ задает непрерывную функцию.

\begin{tm}
Сжимающий оператор подобия $G$ задает непрерывную функцию тогда и только тогда, когда выполнены
следующие условия:
\begin{equation}\label{eq:condD}
\max_{1\leqslant k\leqslant n}|d_k|<1,
\end{equation}
\begin{gather}\label{eq:contf2_b}
\beta_1=f(0)(1-d_1),\\
\beta_k=\sum_{j=1}^{k-1} c_j+f(1)\sum_{j=1}^{k-1}
d_j+f(0)(1-\sum_{j=1}^{k}d_k),\quad k=2,3,\ldots,n,\\
\sum_{j=1}^{n}c_j+(f(1)-f(0))\sum_{j=1}^{n}d_j=f(1)-f(0).\label{eq:contf2_e}
\end{gather}
\end{tm}
\begin{proof}
Оператор подобия $G$ будет сжимающим в пространстве $C[0,1]$ при выполнении условия
$$
|d_k|<1, \quad k=1,2,\ldots, n,
$$
что может быть получено аналогично рассуждениям в лемме~\ref{lem2:1} для пространства
$L_\infty[0,1]$.

Пусть дано разбиение $\{\alpha_k\}_{k=1}^n$. Если функция $f$ является неподвижной точкой
отображения $G$, то $f(\alpha_k)=G_{k+1}(f(0))=G_k(f(1))$, $k=2,3\ldots,n-1$. Кроме того,
$G_1(f(0))=f(0)$ и $G_n(f(1))=f(1)$. Отсюда следует, что самоподобная непрерывная функция $f$
удовлетворяет в точках $\alpha_k$ условиям
\begin{gather}\label{eq:contf1_b}
d_1f(0)+\beta_1=f(0)\quad \text{ в точке } 0=\alpha_1,\\
c_k+d_kf(1)+\beta_k=d_{k+1}f(0)+\beta_{k+1}=f(\alpha_k) \quad \text{ в точках } \alpha_k\quad k=2,3,\ldots,n-1,\\
c_n+d_nf(1)+\beta_n=f(1)\quad \text{ в точке } \alpha_{n-1}.\label{eq:contf1_e}
\end{gather}

Эти соотношения равносильны следующим равенствам
\begin{gather*}
\beta_1=f(0)(1-d_1),\\
\beta_k=\sum_{j=1}^{k-1} c_j+f(1)\sum_{j=1}^{k-1}
d_j+f(0)(1-\sum_{j=1}^{k}d_k),\quad k=2,3,\ldots,n,\\
\sum_{j=1}^{n}c_j+(f(1)-f(0))\sum_{j=1}^{n}d_j=f(1)-f(0).
\end{gather*}

Условия~\eqref{eq:condD} и~\eqref{eq:contf2_b}--\eqref{eq:contf2_e} являются необходимыми и
достаточными условиями, чтобы функция $G(f)$ была непрерывна в точках $\alpha_k$, если исходная
функция $f$ непрерывна. В интервалах $(\alpha_{k-1},\alpha_k)$ непрерывная функция $f$ под
действием сжимающего отображения отображения $G$ перейдет в непрерывную функцию $G(f)$. В
результате при выполнении условий~\eqref{eq:condD}--\eqref{eq:contf2_b} произвольная непрерывная на
отрезке $[0,1]$ функция $f$ перейдет в непрерывную функцию $G(f)$. Таким образом, для произвольной
непрерывной начальной функции $f_0$ последовательность непрерывных функций $f_n=G^{n-1}f_0$
равномерно сходится к неподвижной функции сжимающего отображения и этот предел является непрерывной
функцией.
\end{proof}

\begin{rim} Условия~\eqref{eq:contf2_b}--\eqref{eq:contf2_e} являются обобщением условия
Барнсли (или перекрёстного условия)~\cite{Bar}.
\end{rim}
Условие~\eqref{eq:condD} было получено в~\cite{Bar}.
Условия~\eqref{eq:contf2_b}--\eqref{eq:contf2_e} не описывались, потому что в этой работе решалась
задача интерполяции. А именно, описать все непрерывные аффино-самоподобные функции, графики которых
проходят через через заданные точки $(x_k,y_k)$ $k=1,2,\ldots, n$. Коэффициенты $\{d_k\}$ задаются
произвольно, но подчиненные условию~\eqref{eq:condD}. Остальные параметры операторов подобия
$\{a_k\}$, $\{c_k\}$, $\{\alpha_k\}$ и $\{\beta_k\}$, $k=1,2,\ldots, n$ выражаются через координаты
узлов интерполяции и коэффициенты $\{d_k\}$. В нашем случае самоподобная функция $f$ проходит через
точки $(\alpha_k,\beta_k)$. Через параметры самоподобия не выражаются значения $f(0)$ и $f(1)$.

При описании самоподобных вероятностных мер бывает более удобным рассматривать функции, заданные на
отрезке $[0,1]$ и принимающие на концах значения $f(0)=0$ и  $f(1)=1$. В этом случае
условия~\eqref{eq:contf2_b}--\eqref{eq:contf2_e} принимают более простой вид:
\begin{gather*}
\beta_1=0,\\
\beta_k=\sum_{j=1}^{k-1} (c_j+d_j),\quad k=2,3,\ldots,n,\\
\sum_{j=1}^{n}(c_j+d_j)=1.
\end{gather*}

\section{Связь самоподобных функций из $L_p[0,1]$ и самоподобных мер}
Напомним вкратце конструкцию самоподобных мер. Общую конструкцию фрактальных и самоподобных мер
 можно найти, например, в~\cite{Hu}. Некоторые частные случаи самоподобных сингулярных мер
 представлены в работе~\cite{SV}. Ограничимся рассмотрением вероятностных мер.

Пусть $S_k$, $k=1,2,\ldots,n$
--- семейство сжимающих отображений единичного отрезка $I_0=[0,1]$ в себя, удовлетворяющие условию
$S_k([0,1])=[\alpha_{k},\alpha_{k+1}]$. Пусть также дан набор положительных чисел
$\{\rho_k\}_{k=1}^n$, подчиняющихся условию $\sum_{k=1}^n
\rho_k=1$.

Тогда~(\cite{Hu}), существует и единственна мера $\mu$, заданная уравнением самоподобия
$$
\mu=\sum_{k=1}^n \rho_k \mu\circ S^{-1}_k.
$$
Такие меры называют самоподобными.

С неубывающей непрерывной слева самоподобной функцией $f$ свяжем меру:
$$\mu_f([\zeta,\xi))=f(\xi)-f(\zeta).
$$

Естественно возникает вопрос, является ли эта мера самоподобной.

В данной связи интерес представляет вопрос о нахождении условий в терминах параметров самоподобия
$\{a_k\}$, $\{c_k\}$, $\{d_k\}$, $\{\beta_k\}$, гарантирующих неубывание самоподобной функции.

Рассмотрим ограниченные самоподобные функции, нормированные условиями $f(0)=0$, $f(1)=1$. В этом
случае $\beta_1=0$. Также удобно положить $\beta_{n+1}=1$.

\begin{tm}\label{tm:neubyv}
Чтобы самоподобная непрерывная слева ограниченная функция $f$ была неубывающей необходимо, чтобы
для всех $k=1,2,\ldots,n$ выполнялось $c_k+d_k\geqslant 0$, $\beta_{k}\leqslant \beta_{k+1}$ и
достаточно, чтобы для всех $k=1,2,\ldots,n$ выполнялось $c_k\geqslant 0$, $d_k\geqslant 0$,
$\beta_{k}\leqslant\beta_{k+1}$.
\end{tm}

\begin{proof}
Действительно, из самоподобия функции $f$ следует
\begin{gather*}
f(\alpha_{k}+0)=G_{k}(f(0))=c_{k}\cdot f(0)+\beta_{k}+d_{k}\cdot f(0)=\beta_{k},\\
f(\alpha_{k+1}+0)=G_{k+1}(f(0))=c_{k+1}\cdot f(0)+\beta_{k+1}+d_{k}\cdot f(0)=\beta_{k+1}.
\end{gather*}
Из этих условий вытекает необходимость условий $\beta_{k}\leqslant \beta_{k+1}$. C другой стороны
справедливы соотношения
$$
f(\alpha_k)=G_{k-1}(f(1))=c_{k-1}\cdot f(1)+\beta_{k-1}+d_{k-1}\cdot f(1)=
c_{k-1}+\beta_{k-1}+d_{k-1}\leqslant\beta_{k},
$$
которые приводят к необходимости условий $c_k+d_k\geqslant 0$, $k=1,2,\ldots,n$.

Для получения достаточных условий заметим, что при $\beta_{k-1}\leqslant \beta_k$ выполнено
$f(\alpha_{k-1})\leqslant f(a_k)$, т.\,е. функция неубывает "<глобально"> (в точках разбиения
$\{\alpha_k\}$).  Установим теперь неубывание функции на интервале $(\alpha_k,\alpha_{k+1})$.
Очевидно, что функция, являющаяся неподвижной точкой сжимающего отображения, неотрицательна при
$c_k\geqslant 0$, $d_k\geqslant 0$, $\beta_{k-1}\leqslant \beta_k$. А для произвольной точки
$x\in(\alpha_k,\alpha_{k+1})$ выполнено
$$
f(x)=c_k\cdot t+d_k\cdot f(t)+\beta_k\geqslant f(\alpha_k+0),
$$ где
$x=a_k\cdot t+\alpha_k$.
\end{proof}

Приведем пример, показывающий, что условия $\beta_{k-1}\leqslant \beta_k$, $c_k+d_k\geqslant 0$ не
являются достаточными. Рассмотрим функцию с параметрами самоподобия $n=3$, $a_1=a_2=a_3=1/3$,
$c_1=0$, $d_1=1/2$, $\alpha_1=\beta_1=0$, $c_2=d>0$, $d_2=-d<0$, $\beta_2=1/2$, $c_3=0$, $d_3=1/2$,
$\beta_3=1/2$, где $0<d<1$ --- произвольный параметр. Эта функция является непрерывной и легко
видеть, что все необходимые условия выполнены. Тогда $f(1/3)=1/2$, а $f(1/3+1/9)=d\cdot(1/3)-d\cdot
f(1/3)+1/2=1/2-d/6<1/2$.

Рассмотрим неубывающие ограниченные самоподобные функции $f$, нормированные условиями $f(0)=0$,
$f(1)=1$. Кроме того, пусть $c_k=0$, $k=1,2,\ldots,n$. В этом случае из теоремы~\ref{tm:neubyv}
следует, что необходимыми и достаточными условиями неубывания функции будут условия $d_k\geqslant
0$, $\beta_k\leqslant\beta_{k+1}$, $k=1,2\ldots,n$. Легко заметить, что мера, построенная по такой
самоподобной функции является самоподобной. При этом роль чисел $\rho_k$ выполняют положительные
числа $d_k$, а отображения $S^{-1}_k\colon[\alpha_{k},\alpha_{k+1}]\to [0,1]$ в данном случае
действуют по правилу
$$
S^{-1}_k(x)=\dfrac{x-\alpha_k}{a_k}.
$$

\begin{rim}
Рассмотрение неубывающих самоподобных функций с ненулевыми параметрами $\{c_k\}$ приводит к
самоподобным мерам с абсолютно-непрерывной частью. Подробное рассмотрение таких мер выходит за
рамки данной статьи и является темой дальнейшего исследования.
\end{rim}

\subsection{Непрерывные самоподобные функции с неограниченной вариацией}

Рассмотрим непрерывные самоподобные функции, заданные некоторым набором самоподобия $\{a_k\}$,
$\{c_k\}$, $\{d_k\}$, $\{\beta_k\}$ и нормированные условием $f(0)=0$, $f(1)=1$. Изучим только
частный случай таких функций, а именно, предположим, что $c_k=0$ $k=1,2,\ldots,n$. В этом случае из
условий~\eqref{eq:condD}--\eqref{eq:contf2_e} следует, что
$$
 \beta_1=0,\quad  \beta_{k+1}-\beta_k=d_k, \quad \sum_{k=1}^n d_k=1,\quad \max_{1\leqslant k\leqslant n}\{|d_k|\}<1.
$$
Введем величину $D=\sum_{k=1}^n |d_k|$. Нетрудно видеть, что из условия $\sum_{k=1}^n d_k=1$
следует, что $D\geqslant 1$.

\begin{tm}
Самоподобная непрерывная функция $f$, с параметрами самоподобия $\{a_k\}$, $\{c_k\}$, $\{d_k\}$,
$\{\beta_k\}$, удовлетворяющая условиям $f(0)=1$, $f(1)=1$, $c_k=0$ $k=1,2,\ldots,n$, имеет
ограниченную вариацию тогда и только тогда, когда $D\leqslant 1$.
\end{tm}
\begin{proof}
\emph{Шаг 1.}
Построим последовательность  $T_m$ разбиений отрезка $[0,1]$. Пусть $T_1=\{\alpha_k\}_{k=1}^{n+1}$
($\alpha_{0}=0,\alpha_{n+1}=1$) --- исходное разбиение отрезка.   $T_m=\{a_k\cdot x+\alpha_k,x\in
T_{m-1}, k=1,2,\ldots,n\}$. Точки разбиения $T_m$ делят отрезок $[0,1]$ на $n^m$ подотрезков.

Каждому отрезку, образованному разбиением $T_m$, сопоставим последовательность чисел
$k_1,k_2,\ldots,k_m$, $k_i=1,2,\ldots,n$, $i=1,2,\ldots,m$ (номер отрезка). Отрезку разбиения $T_1$
соответствует номер $k_1$, если границы этого отрезка строятся по правилу
$$
x_{k_1}=\alpha_{k_1}, \quad x_{k_1+1}=\alpha_{k_1+1}, k_1=1,2,\ldots,n.
$$
Отрезку разбиения $T_2$ соответствует пара $k_1,k_2$, если его концы $x_{k_2}$ и $x_{k_2+1}$
построены по правилу
$$
x_{k_2}=a_{k_2}x_{k_1}+\alpha_{k_2}, \quad x_{k_2+1}=a_{k_2}x_{k_1+1}+\alpha_{k_2},
k_1,k_2=1,2,\ldots,n,
$$
где $[x_{k_1},x_{k_1+1}]$ --- отрезок разбиения $T_1$ с номером $k_1$. Далее построение
последовательности номеров отрезка проводится по индукции: если каждому отрезку разбиения $T_{m-1}$
соответствует последовательность $k_1,k_2,\ldots,k_{m-1}$, то отрезок разбиения $T_m$ имеет номер
$k_1,k_2,\ldots,k_m$, если его концы $x_{k_m}$ и $x_{k_m+1}$ построены по правилу
$$
x_{k_m}=a_{k_m}x_{k_{m-1}}+\alpha_{k_m}, \quad x_{k_m+1}=a_{k_m}x_{k_{m-1}+1}+\alpha_{k_m},
$$
где $[x_{k_{m-1}},x_{k_{m-1}+1}]$ --- отрезок с номером $k_1,k_2,\ldots,k_{m-1}$,
$k_i=1,2,\ldots,n$, $i=1,2,\ldots,m-1$. Нетрудно заметить, что длина отрезка с номером
$k_1,k_2,\ldots,k_m$ равна $a_{k_1}a_{k_2}\ldots a_{k_m}$.

\emph{Шаг 2.}
Пусть $f_0(x)=x$ на $[0,1]$. Последовательность $f_m=G(f_{m-1})$ сходится к предельной самоподобной
функции. Докажем, что на отрезке с номером $k_1,k_2,\ldots,k_m$ функция $f_m$ линейна и
\begin{equation}\label{eq:fm}
f_m(x)=\dfrac{d_{k_1}d_{k_2}\ldots d_{k_m}}{a_{k_1}a_{k_2}\ldots a_{k_m}}(x-x_{k_m})+
\sum_{j=1}^m\beta_{k_j}\prod_{i=j+1}d_{k_i}.
\end{equation}

Очевидно, что по определению действия оператора подобия $G$ на непрерывные
функции~\eqref{eq:nepr_otobr}, кусочно-линейная функция перейдет в кусочно-линейную же функцию.
Пусть $[x_{k_1},x_{k_1+1}]$ --- отрезок разбиения $T_1$ с номером $k_1$. В силу того же
свойства~\eqref{eq:nepr_otobr} имеем
\begin{gather*}
f_1(x_{k_1})=d_{k_1}f_0(0)+\beta_{k_1}=\beta_{k_1},\\
f_1(x_{k_1+1})=d_{k_1}f_0(1)+\beta_{k_1}=\beta_{k_1+1}.
\end{gather*}
Таким образом, на отрезке разбиения $T_1$ с номером $k_1$ функция $f_1$ линейна и на его концах
принимает значения $\beta_{k_1}$ и $\beta_{k_1+1}$, т.\,е. определяется формулой
$$
f_1(x)=\dfrac{d_{k_1}}{a_{k_1}}(x-x_{k_1})+\beta_{k_1}.
$$

Дальнейшее рассуждение проведем по индукции. Пусть на отрезке разбиения $T_m$ с номером
$k_1,k_2,\ldots,k_{m}$ функция $f_m$ имеет вид~\eqref{eq:fm}. Концы отрезка разбиения $T_{m+1}$ с
номером $k_1,k_2,\ldots,k_{m},k_{m+1}$ строятся по правилу
$$
x_{k_{m+1}}=a_{k_{m+1}}x_{k_m}+\alpha_{k_{m+1}},\quad
x_{k_{m+1}+1}=a_{k_{m+1}}x_{k_m+1}+\alpha_{k_{m+1}},
$$
где $[x_{k_m},x_{k_m+1}]$ --- отрезок с номером $k_1,k_2,\ldots,k_{m}$. На этом отрезке $f_{m+1}$
кусочно-линейна и на его концах принимает значения
\begin{gather*}
f_{m+1}(x_{k_{m+1}})=d_{k_{m+1}}f_m(x_{k_m})+\beta_{k_{m+1}}=\sum_{j=1}^m\beta_{k_j}\prod_{i=j+1}d_{k_i}
+\beta_{k_{m+1}}\\
f_{m+1}(x_{k_{m+1}+1})=d_{k_{m+1}}f_m(x_{k_m+1})+\beta_{k_{m+1}}=d_{k_1}d_{k_2}\ldots
d_{k_m}d_{k_{m+1}}+\sum_{j=1}^m\beta_{k_j}\prod_{i=j+1}d_{k_i} +\beta_{k_{m+1}},
\end{gather*}
откуда и следует справедливость утверждения.

\emph{Шаг 3.} Заметим, что в точках $x_{k_1}$ разбиения $T_1$ выполнено
\begin{equation}\label{eq:sovpt1}
f_1(x_{k_1})=d_{k_1}f_0(0)+\beta_{k_1}=\beta_{k_1},\quad f(x_{k_1})=d_{k_1}f(0)+\beta_{k_1},
\end{equation}
т.\,е.
$$
f_1|_{T_1}=f|_{T_1}.
$$

А для точек $x_{k_m}$ разбиения $T_m$ выполнено
$$
f_m(x_{k_m})=d_{k_m}f_{m-1}(x_{k_{m-1}})+\beta_{k_m},\quad
f(x_{k_m})=d_{k_m}f(x_{k_{m-1}})+\beta_{k_m},
$$
где ${x_{k_m}}$ и ${x_{k_{m-1}}}$ --- левые концы отрезков с номерами $k_1,k_2,\ldots,k_{m}$ и
$k_1,k_2,\ldots,k_{m-1}$ соответственно. Из этих равенств с учетом~\eqref{eq:sovpt1} следует, что
в точках разбиения $T_m$ выполнены равенства
$$
f_m|_{T_m}=f|_{T_m}.
$$
Следовательно, вариации функций $f$ и $f_m$ по разбиению $T_m$ совпадают.

Вариация функции $f_m$ по разбиению $T_m$ равна
$$
\Var_{T_m}f_m=\sum_{1\leqslant k_1,k_2,\ldots, k_m\leqslant n }
|d_{k_1}d_{k_2}\ldots d_{k_m}|=D^m,
$$

Так как (см.~\cite{Nat}, теорема 2, стр.~211) $\Var_{T_m}f\to \Var_0^1 f$ при $m\to\infty$, то
$\Var_{T_m}f_m\to \Var_0^1 f$, где $\Var_{T_m}f$ это вариация функции $f$ по разбиению $T_m$.

Следовательно,
$$
\Var_0^1 f=\lim_{m\to\infty} D^m.
$$
\end{proof}

Таким образом, в классе непрерывных самоподобных функций, удовлетворяющих условиям $f(0)=0$ и
$f(1)=1$, а параметры самоподобия которой подчинены условиям $c_k=0$, $k=1,2,\ldots,n$,
ограниченную вариацию имеют только те функции, для которых $D=\sum_{k=1}^n |d_k|=1$.

\begin{rim}
Для фрактальных кривых, заданных аффинными (не обязательно сжимающими) операторами, условию
$D\leqslant 1$ соответствует условие $\rho_1\leqslant 1$~\cite{Prot1}.
\end{rim}

\section{Примеры}\label{pt:4}
\subsection{Некоторые конкретные самоподобные функции}

1) Характеристическая функция интервала $(\zeta,\xi)\subset [0,1]$ является самоподобной функцией с
параметрами самоподобия $n=3$, $a_1=\zeta$, $a_2=\xi-\zeta$, $a_3=1-\xi$,
$c_1=c_2=c_3=d_1=d_2=d_3=\beta_1=\beta_3=0$, $\beta_2=1$.

2) Кусочно-постоянная функция $f(x)=\sum_{k=1}^n s_k\chi_{[\alpha_{k},\alpha_{k+1})}(x)$, где
$0=\alpha_1<\alpha_2<\ldots<\alpha_n<\alpha_{n+1}=1$ --- разбиение отрезка $[0,1]$, также является
самоподобной функцией с параметрами самоподобия $c_1=c_2=\ldots=c_n=0$, $d_1=d_2=\ldots=d_n=0$,
$a_k=\alpha_{k+1}-\alpha_k$, $\beta_k=s_k$, $k=1,2,\ldots,n$.

3) Рассмотрим двупараметрическое семейство непрерывных функций  $f_{a,\delta}$, где $a\in(0,1/2)$ и
$\delta\in[0,1/3)$, определяемое параметрами самоподобия
\begin{multline*}
n=3,\quad a_1=a_2=a,\quad a_2=1-2a,\\ d_1=d_3=1/2+\delta,\quad d_2=-2\delta,
\quad\beta_1=0,\quad\beta_2=d_1=1/2+\delta,\quad\beta_3=d_1+d_2=1-\delta.
\end{multline*}
В частности, $f_{1/3,0}$ представляет собой хорошо известную функцию Кантора. Для этих функций
$D=\sum_{k=1}^3|d_k|=1+4\delta$ и при $\delta>0$ выполнено $D>1$, т.\,е. функции имеют
неограниченные вариации.

Автор выражает благодарность В.~Ю.~Протасову за ряд ценных замечаний и полезные советы.


\vspace{0.3cm}




\begin{thebibliography}{99}
\bibitem{Hu} J.~Hutchinson, \emph{Fractals and
Self-similatity}//Indiana University Math. J., \textbf{30} (1981), pp.~713--741.

\bibitem{Bar} M.~Barnsley, \emph{Fractals everywere}//Academic Press, 1988.

\bibitem{SV} M.~Solomyak, E.~Verbitsky, \textit{On a spectral problem
related to self-similar measures}//Bull. London Math. Soc.,
\textbf{27}\,(1995), pp.~242--248.

\bibitem{Str} R.S.Strihartz, \emph{Self-similar measures and their Fourier transform,
I}// Indiana University Math.J.,\textbf{39} (1990) pp.797--817.

\bibitem{Zig} А.Зигмунд, \emph{Тригонометрические ряды, том 1}, М., Мир,1965.

\bibitem{deRam1} De Rham G., \emph{Une peu de math\'ematique \`a propos d'une courbe plane}// Rev. de
math. elemetaires, \textbf{2}:4,5 (1947), pp.73--76, 89--97.

\bibitem{deRam2} De Rham G., \emph{Sur une courbe plane}// J. Math. Pures Appl. (9), \textbf{35}
(1956), pp. 25--42.

\bibitem{deRam3} De Rham G., \emph{Sur les courbes limit\'es de polygones obtenus par
trisection}//Enseign. Math., (2), \textbf{5} (1959), pp.29--43.

\bibitem{CavDah} Cavaretta A.S., Dahmen W., Micchelli C.A., \emph{Stationary
dubdivision}//Mem. Amer. Math. Soc., \textbf{93}:453 (1991), p.186.

\bibitem{Derf} Г.~А.~Дерфель \emph{Вероятностные методы для одного класса функционально-разностных
уравнений}//Укр. матем. журн., \textbf{41}:10 (1989), стр. 1137--1141.

\bibitem{Nik} П.~П.~Никитин, \emph{Хаусдорфова размерность гармонической меры на кривой де
Рама}//Зап. науч. семинара ПОМИ, \textbf{283} (2001), стр. 206--223.

\bibitem{Daub1} Daubechies I.,Lagarias J., \emph{Two scale difference equations. I. Existence and global regularity of
solutions}//SIAM. J. Anal., \textbf{22}:5 (1991), pp. 1388--1410.

\bibitem{Daub2} Daubechies I.,Lagarias J., \emph{Two scale difference equations. I. Local regularity,
infinite products of matrices and fractals}//SIAM. J. Anal., \textbf{23}:4 (1992), pp. 1031--1079.

\bibitem{Prot1} Protasov V. \emph{Refinement equations with nonnegative coefficients}//J. Fourier
Anal. Appl., \textbf{6}:1, (2000), pp. 55--78.

\bibitem{Prot2} В.~Ю.~Протасов \emph{Фрактальные кривые и всплески}//Изв.РАН. Серия матем.,
\textbf{70}:5 (2006), стр.105--145.

\bibitem{LauWang} Lau K.-S., Wang J., \emph{Characterization of $L_p$-solutions for two-scale dilation
equations}//SIAM. J. Math. Anal., \textbf{26}:4, (1995), 1018--1046.

\bibitem{VlSh} А.~А.~Владимиров, И.~А.~Шейпак, \emph{Самоподобные функции в пространстве \(L_2[0,1]\) и задача
Штурма-Лиувилля с сингулярным весом}//http://www.arxiv.org/math.FA/0405410

\bibitem{KacKr} И.~С.~Кац, М.~Г.~Крейн., \textit{О спектральных функциях
струны.} В кн.:~Ф.~Аткинсон. \textit{Дискретные и непрерывные граничные задачи.} М., 1968,
стр.~648--733.

\bibitem{Nat} И.~П.~Натансон. \emph{Теория функций вещественной переменной}//М.:\,Наука, 1974.


\end{thebibliography}
\end{document}